\documentclass[12pt]{article}
\usepackage{amsfonts}
\usepackage{latexsym,amsmath}
\usepackage{amssymb,array}
\makeatletter \oddsidemargin 0in \evensidemargin 0in \textwidth
16cm \RequirePackage[dvips]{graphicx} \textheight 20cm
\newtheorem{t1}{Theorem}[section]
\newtheorem{p1}{Proposition}[section]

\newtheorem{c1}{Corollary}[section]
\newtheorem{d1}{Definition}[section]
\newtheorem{r1}{Remark}[section]

\newtheorem{ex}{Example}[section]

\begin{document}
\title{Bivariate Extension of (Dynamic) Cumulative Past Entropy}
\author{{\bf Amarjit Kundu}\\
Department of Mathematics\\Santipur College\\P.O. Santipur, Dist. Nadia, W.B., India
\and
{\bf Chanchal Kundu}\footnote{Corresponding
author e-mail: ckundu@rgipt.ac.in, chanchal$_{-}$kundu@yahoo.com.}\\
Department of Mathematics\\Rajiv Gandhi Institute of Petroleum Technology\\Rae Bareli 229 316, U.P., India}
\date{Revised version to appear in {\it Communications in Statistics$-$Theory \& Methods}, $\copyright$ by Taylor $\&$ Francis Group.\\Submitted: August, 2014}
\maketitle
\begin{abstract} Recently, the concept of cumulative residual entropy (CRE) has been studied by many researchers in higher dimensions. In this article, we extend the definition of (dynamic) cumulative past entropy (DCPE), a dual measure of  (dynamic) CRE, to bivariate setup and obtain some of its properties including bounds. We also look into the problem of extending DCPE for conditionally specified models. Several properties, including monotonicity, and bounds of DCPE are obtained for conditional distributions. It is shown that the proposed measure uniquely determines the distribution function. Moreover, we also propose a stochastic order based on this measure.
\end{abstract}
{\bf Key Words and Phrases:} Cumulative past entropy, bivariate reversed hazard rate and expected inactivity time, stochastic ordering.\\
{\bf AMS 2010 Classifications:} Primary 62G30; Secondary 62E10, 62B10, 94A17.
\section{Introduction}
In recent years, cumulative residual entropy and cumulative past entropy have been considered as a new measure of information that parallel Shannon (1948) entropy. Let $X$ be an absolutely continuous nonnegative random variable with survival function $\overline F(\cdot)=1-F(\cdot)$ and probability density function (pdf) $f(\cdot)$. Then Shannon's differential entropy is defined as
\begin{eqnarray}\label{eq1.1}
H(X)=-\int_{0}^\infty f(x)\log f(x) dx,
\end{eqnarray}
where `log' means natural logarithm and, by convention, $0\log0=0$. It measures the expected uncertainty contained in $f$ about the predictability of an outcome of $X$. In spite of enormous success of Shannon entropy, the differential entropy (\ref{eq1.1}) presents various deficiencies when it is used as a continuous counterpart of the classical Shannon entropy for discrete random variable. Several attempts have been made in order to define possible alternative information measures. Recently, Rao et al. (2004) identified some limitations of the use of (\ref{eq1.1}) in measuring randomness of certain systems and introduced an alternative measure of uncertainty called cumulative residual entropy (CRE) defined as
\begin{eqnarray}\label{eq1.2}
\varepsilon(X)=-\int_{0}^{\infty}\overline F(x)\log{\overline
F(x)}dx,
\end{eqnarray}
which relates to uncertainty on the future lifetime of a system. This measure is based on survival function and is particularly suitable to describe the information in problems related to ageing properties in the reliability theory. Motivated by the salient features of (\ref{eq1.2}), Di Crescenzo and Longobardi (2009) proposed a dual concept of CRE called cumulative past entropy (CPE) defined as
\begin{eqnarray}\label{eq1.3}
\overline\varepsilon(X)=-\int_{0}^{\infty}F(x)\log{F(x)}dx,
\end{eqnarray}
which measures information concerning past lifetime.\\
\hspace*{.2in} Length of time during a study period has been considered as a prime variable of interest in many areas such as reliability, survival analysis, economics, business, etc. In particular, consider an item under study, then the information about the residual (past) lifetime is an important task in many applications. In such cases the information measures are functions of time, and thus they are dynamic in nature. Asadi and Zohrevand (2007) further studied the function obtained from (\ref{eq1.2}) in the residual setup called dynamic CRE (DCRE), given by
\begin{eqnarray}\label{eq1.4}
\varepsilon_X(t)=-\int_{t}^{\infty}\frac{\overline F(x)}{\overline F(t)}\log\frac{\overline F(x)}{\overline F(t)}dx.
\end{eqnarray}
For a discussion on the properties and generalization of (dynamic) CRE one may refer to Rao (2005), Navarro et al. (2010), Abbasnejad et al. (2010), Kumar and Taneja (2011), Navarro et al. (2011), Sunoj and Linu (2012), Khorashadizadeh et al. (2013), Psarrakos and Navarro (2013), Navarro et al. (2014), Chamany and Baratpour (2014), among others. In analogy with (\ref{eq1.4}), Di Crescenzo and Longobardi (2009) also studied CPE for past lifetime called dynamic CPE (DCPE), defined as
\begin{eqnarray}\label{eq1.5}
\overline\varepsilon_X(t)=-\int_{0}^{t}\frac{F(x)}{F(t)}\log\frac{F(x)}{F(t)}dx.
\end{eqnarray}
For more properties, applications and recent developments of (\ref{eq1.5}), one may refer to Abbasnejad (2011) and Di Crescenzo and Longobardi (2013).

\hspace*{.2in} In studying the reliability aspects of multi-component system with each component having a lifetime depending on the lifetimes of the other components, multivariate life distributions are employed. Reliability characteristics in the univariate case can be extended to the corresponding multivariate version. In a recent work, Rajesh et al. (2014) have considered extension of DCRE to bivariate setup and study its properties. Some bivariate distributions are also characterized there. Several generalizations to the concept of bivariate DCRE can be found in Sunoj and Linu (2012) and Rajesh et al. (2014). In various contexts the uncertainty is not necessarily related to the future but may refer to the past. Even though a lot of interest has been evoked on the bivariate extension of CRE, no works for CPE, to the best of our knowledge, till now, seem to have been done in higher dimension. It is to be noted that the concepts in past time are more appropriate than those truncated from below when the observations are predominantly from left tail. This shows the relevance and usefulness of studying CPE when uncertainty is related to the past.\\
\hspace*{.2in} In the present paper we consider (dynamic) CPE for bivariate setup, and study its various properties useful in reliability modeling. The rest of the paper is arranged as follows. Section 2 includes the definition and basic properties of bivariate CPE. Some bounds for bivariate CPE are also obtained. In Section 3 we look into the behavior of dynamic CPE for conditional distributions. Several properties of the measures are studied along with monotonicity and some characterization theorems arising out of it. A stochastic order is proposed and studied which is based on the measures.
\section{Definition and properties of bivariate CPE}
The topic of measuring the information content for bivariate (multivariate) distributions when their supports are truncated progressively are considered in recent past. A significant results in this area have been provided in Ebrahimi et al. (2007). In this section we look into the problem of extending CPE defined in (\ref{eq1.3}) to the bivariate setup. Let $X=(X_1,X_2)$ be a random vector of nonnegative absolutely continuous random variables with joint pdf $f(x_1,x_2)$. We may think of $X_i,~i=1,2$, as the lifetimes of the members of a group or components of a system. Then the Shannon differential entropy of $(X_1,X_2)$ is defined as
\begin{eqnarray}\label{eq2.0}
H(X_1,X_2)=-\int_0^\infty dx_1\int_0^\infty f(x_1,x_2)\log f(x_1,x_2)dx_2.
\end{eqnarray}
One of the main problems encountered while extending a univariate concept to higher dimensions is that it cannot be done in a unique way. A natural extension of CPE (\ref{eq1.3}) to the bivariate setup can be obtained from (\ref{eq2.0}) by replacing $f(x_1,x_2)$ by $F(x_1,x_2)$ as given in the following definition. Since the past lifetime has always a finite support we restrict our attention to random variables with finite supports. Therefore, we assume that the support of $(X_1,X_2)$ is included in $(0,b_1)\times(0,b_2)$ for some nonnegative real values $b_1,b_2$.
\begin{d1} Let $X=(X_1,X_2)$ be a nonnegative bivariate random vector admitting an absolutely continuous distribution function with joint pdf $f(x_1,x_2)$, distribution function $F(x_1,x_2)$, marginal distribution functions $F_i(x_i)$ and marginal pdfs $f_{X_i}(x_i),~i=1,2$. We define the bivariate CPE as
\begin{eqnarray}\label{eq2.1}
\overline\varepsilon(X_1,X_2)=-\int_0^{b_1}\int_0^{b_2} F(x_1,x_2)\log F(x_1,x_2)dx_2dx_1
\end{eqnarray}
provided the integral on the right hand side is finite.
\end{d1}
\hspace*{.2in} The following example clarifies the effectiveness of the proposed measure.
\begin{ex}\label{ex2.1} Let $X_1$ and $X_2$ are random lifetimes of two components with joint pdf
\begin{equation*}f(x_1,x_2)=\left\{\begin{array}{ll}2,~0<x_1<1,~0<x_2<x_1\\
0,~{\rm otherwise}.
\end{array}\right.
\end{equation*}
Then, $H(X_1,X_2)=-\log2,~H(X_1)=H(X_2)=1/2-\log2$ and $\overline\varepsilon(X_1,X_2)=(1-\log2)/4,~\overline\varepsilon(X_1)=2/9,~\overline\varepsilon(X_2)=\frac{10}{9}-\frac{4}{3}\log2.$ Here $H(X_1,X_2)$ along with $H(X_1)$ and $H(X_2)$ are negative which do not make a sense whereas their CPE measures are positive. Most importantly, $H(X_1)$ and $H(X_2)$ are identical but $\overline\varepsilon(X_1)$ and $\overline\varepsilon(X_2)$ are different. $\hfill\square$
\end{ex}
Indeed, since $\log x\leqslant x-1$ for all $x>0$, we have
$$\overline\varepsilon(X_1,X_2)\geqslant\int_0^{b_1}\int_0^{b_2} F(x_1,x_2)\left(1-F(x_1,x_2)\right)dx_2dx_1,$$
where $1-F(x_1,x_2)={\rm P}(X_1>x_1~{\rm or}~X_2>x_2)$, which gives the probability that at least one of the components will survive beyond the time $(x_1,x_2)$.
In connection with Example \ref{ex2.1} the above right hand side expression is evaluated as 0.028 which verifies the proposed lower bound of bivariate CPE.\\
We also recall from (\ref{eq2.1}) that if $X_1$ and $X_2$ are independent then
\begin{eqnarray}\label{eq2.2}
\overline\varepsilon(X_1,X_2)=\left[\int_0^{b_2} F_2(x_2)dx_2\right]\overline\varepsilon(X_1)+\left[\int_0^{b_1} F_1(x_1)dx_1\right]\overline\varepsilon(X_2).
\end{eqnarray}
The following additive property of bivariate CPE is due to Di Crescenzo and Longobardi (2009).
\begin{p1}Let $X=(X_1,X_2)$ be a nonnegative bivariate random vector where $X_1,X_2$ are independent random variables with supports $[0,b_1]$ and $[0,b_2]$, respectively. Then
$$\overline\varepsilon(X_1,X_2)=\left(b_2-\mu_2\right)\overline\varepsilon(X_1)+\left(b_1-\mu_1\right)\overline\varepsilon(X_2),$$
where $\mu_i=E(X_i)$ and $b_i$ are finite. In particular, if $X_1$ and $X_2$ both have support $[0,b]$ and expectation $\mu$, then we have
$$\overline\varepsilon(X_1,X_2)=(b-\mu)\left[\overline\varepsilon(X_1)+\overline\varepsilon(X_2)\right]$$
which shows that bivariate CPE also have an appealing property in analogy with Shannon's differential entropy for two-dimensional random variable e.g., $H(X_1,X_2)=H(X_1)+H(X_2)$, if $X_1$ and $X_2$ are independent.
\end{p1}

\hspace*{.2in} Let us analyze the effect of linear transformations on the bivariate CPE. The proof is immediate from (\ref{eq2.1}).
\begin{p1} Let $Y=(Y_1,Y_2)$ be a nonnegative bivariate random vector where $Y_i=c_iX_i+d_i$ with $c_i>0$ and $d_i\geqslant0$. Then
$$\overline\varepsilon(Y_1,Y_2)=c_1c_2\overline\varepsilon(X_1,X_2),$$
which shows that bivariate CPE is a shift-independent measure.
\end{p1}
Now we show that bivariate CPE is not invariant under non-singular transformations.
\begin{p1} Let $Y=(Y_1,Y_2)$ be a nonnegative bivariate random vector. If $Y_i=\varphi_i(X_i),~i=1,2$ are one-to-one transformations with $\varphi_i(x_i)$ are differentiable functions, then
$$\overline\varepsilon(Y_1,Y_2)=-\int_0^{b_1}\int_0^{b_2} F(x_1,x_2)\log F(x_1,x_2)|J|dx_2dx_1,$$
where $J=\frac{\partial}{\partial x_1}\varphi_1(x_1)\frac{\partial}{\partial x_2}\varphi_2(x_2)$ is the Jacobian of the transformation.
\end{p1}

\hspace*{.2in} A sharper lower bound for bivariate CPE is given in the following theorem.
\begin{t1} For a nonnegative bivariate random vector $X=(X_1,X_2)$
$$\overline\varepsilon(X_1,X_2)\geqslant {\rm max}\left[C_1e^{H(X_1)},C_2e^{H(X_2)}\right],$$
where $$C_1=\exp\left[E_{X_1}\left(\log\int_0^{b_2} F(X_1,x_2)|\log F(X_1,x_2)|dx_2\right)\right]$$
$${\rm and}~~ C_2=\exp\left[E_{X_2}\left(\log\int_0^{b_1} F(x_1,X_2)|\log F(x_1,X_2)|dx_1\right)\right].$$
\end{t1}
Proof: Using log-sum inequality we get
\begin{eqnarray*}
&\int_0^{b_1}\left[f_{X_1}(x_1)\log\frac{f_{X_1}(x_1)}{\int_0^{b_2} F(x_1,x_2)|\log F(x_1,x_2)|dx_2}\right]dx_1\\&\geqslant\int_0^{b_1} f_{X_1}(x_1)dx_1\log\frac{\int_0^{b_1} f_{X_1}(x_1)dx_1}{\int_0^{b_1}\int_0^{b_2} F(x_1,x_2)|\log F(x_1,x_2)|dx_2dx_1}=-\log\overline\varepsilon(X_1,X_2),
\end{eqnarray*}
which on simplification reduces to
$$\overline\varepsilon(X_1,X_2)\geqslant C_1e^{H(X_1)},$$
where $C_1=\exp\left[E_{X_1}\left(\log\int_0^{b_2} F(X_1,x_2)|\log F(X_1,x_2)|dx_{2}\right)\right]$. Proceeding analogously one can obtain the required result. $\hfill\square$\\

\hspace*{.2in} If $X=(X_1,X_2)$ represents the lifetimes of two components in a system where both the components are found failed at times $t_1$ and $t_2$, respectively, then, the measure of uncertainty associated with the past lifetimes of the system, called bivariate dynamic CPE, is given by
\begin{eqnarray}\label{eq2.4}\overline\varepsilon_X(t_1,t_2)=-\int_0^{t_1}\int_0^{t_2}\frac{F(x_1,x_2)}{F(t_1,t_2)}\log\frac{F(x_1,x_2)}{F(t_1,t_2)}dx_2dx_1,
\end{eqnarray}
which can be thought of as two-dimensional extension of dynamic CPE. If $X_1$ and $X_2$ are independent, then $$\overline\varepsilon_X(t_1,t_2)=m_2(t_2)\overline\varepsilon_{X_1}(t_1)+m_1(t_1)\overline\varepsilon_{X_2}(t_2)$$ where $m_i(t_i)=E\left(t_i-X_i|X_i\leqslant t_i\right)$, are the marginal expected inactivity time (EIT) of the components $X_i,~i=1,2$. The bivariate dynamic CPE is also not invariant under non-singular transformations.
\begin{p1}Let $Y=(Y_1,Y_2)$ be a nonnegative bivariate random vector. If $Y_i=\varphi_i(X_i),~i=1,2$ are one-to-one transformations with $\varphi_i(x_i)$ are differentiable functions, then
$$\overline\varepsilon_Y(\varphi_1(t_1),\varphi_2(t_2))=-\int_0^{t_1}\int_0^{t_2}\frac{F(x_1,x_2)}{F(t_1,t_2)}\log\frac{F(x_1,x_2)}{F(t_1,t_2)}|J|dx_2dx_1.$$
In particular, if we choose $\varphi_i(X_i)=c_iX_i+d_i$ with $c_i>0$ and $d_i\geqslant0$ for $i=1,2$ then $\overline\varepsilon_Y(\varphi_1(t_1),\varphi_2(t_2))=c_1c_2\overline\varepsilon_X(t_1,t_2)$.
\end{p1}
\hspace*{.2in} When we consider bivariate measures, it is necessary that the measurement on the basis of one component is not affected by the missing or unreliable data on the other component and hence it is necessary to consider component-wise CPE subject to the condition that both the components are found failed at some specified times. Such a measure will be more reliable as the unreliable data is omitted. With this motivation, we now look into the behavior of dynamic CPE for conditional distributions.
\section{Conditional dynamic CPE}
Specification of the joint distribution through its component densities, namely marginals and conditionals has been a problem dealt with by many researchers in the past. It is well known that in general, the marginal densities cannot determine the joint density uniquely unless the variables are independent. Apart from the marginal distribution of $X_i$ and the conditional distribution of $X_j$ given $X_i=t_i$, $i=1,2,~i\neq j$, from which the joint distribution can always be found, the other quantities that are of relevance to the problem are (a) the two conditional distributions of $X_i$ given that $X_j<t_j$, $i,j=1,2,~i\neq j$, (b) marginal and conditional distributions of the same component viz. $X_1$ and the $X_1$ given $X_2=t_2$ or $X_2$ and the $X_2$ given $X_1=t_1$. Characterization of the bivariate density given the forms of the marginal density of $X_1~(X_2)$ and the conditional density of $X_1$ given $X_2=t_2~(X_2$ given $X_1=t_1)$ for certain classes of distributions, have been considered by Seshadri and Patil (1964), Nair and Nair (1988), Hitha and Nair (1991), Arnold et al. (2001) and Navarro and Sarabia (2013). Accordingly in the following Subsections 3.1 and 3.2, we consider conditional dynamic CPE of $X_i$ given $X_j<t_j$ and $X_i$ given $X_j=t_j,~i,j=1,2,~ i\neq j$, respectively and study some characteristics relationships in the context of reliability modeling.
\subsection{Conditional dynamic CPE for $X_i$ given $X_j<t_j$}
Here we look into the behavior of dynamic CPE for conditional distributions. Consider the random variables $\widetilde{Y_i}=(X_i|X_1<t_1,X_2<t_2),~i=1,2,$ which correspond to the conditional distributions of $X_i$ subject to the condition that failure of the first component had occurred in $(0,t_1)$ and the second has failed before $t_2$. The distribution functions of $\widetilde{Y_i}$ are given by $P\left(\widetilde{Y_1}\leqslant y_1\right)=\frac{F(y_1,t_2)}{F(t_1,t_2)},~0<y_1<t_1$ and $P\left(\widetilde{Y_2}\leqslant y_2\right)=\frac{F(t_1,y_2)}{F(t_1,t_2)},~0<y_2<t_2$. Then the dynamic CPE for $\widetilde{Y_i}$, called conditional dynamic CPE (CDCPE), takes the form
\begin{eqnarray}\label{eq3.1}
\overline\varepsilon^*_1(X;t_1,t_2)=-\int_0^{t_1}\frac{F(x_1,t_2)}{F(t_1,t_2)}\log\frac{F(x_1,t_2)}{F(t_1,t_2)}dx_1
\end{eqnarray}
\begin{eqnarray}\label{eq3.2}
{\rm and}~~\overline\varepsilon^*_2(X;t_1,t_2)=-\int_0^{t_2}\frac{F(t_1,x_2)}{F(t_1,t_2)}\log\frac{F(t_1,x_2)}{F(t_1,t_2)}dx_2.
\end{eqnarray}
If $X=(X_1,X_2)$ represents a bivariate random vector, recalling (\ref{eq1.3}), then $\overline\varepsilon^*_1(X;t_1,t_2)$ identifies with the CPE of $(X_1|X_1<t_1,X_2<t_2)$, with a similar interpretation for $\overline\varepsilon^*_2(X;t_1,t_2)$. In particular, if $X_1$ and $X_2$ are independent, then $\overline\varepsilon^*_i(X;t_1,t_2)=\overline\varepsilon_{X_i}(t_i),~i=1,2$. In the sequel we give the definitions of bivariate reversed hazard rate and bivariate EIT functions. For more details one may refer to Roy (2002) and Nair and Asha (2008).
\begin{d1}\label{def3.1} For a random vector $X=(X_1,X_2)$ with distribution functions $F\left(t_1,t_2\right)$
\item[(i)] the bivariate reversed hazard rate is defined as a vector, $\phi^X(t_1,t_2)=\left(\phi_1^X(t_1,t_2),\phi_2^X(t_1,t_2)\right)$ where $\phi_i^X(t_1,t_2)=\frac{\partial}{\partial t_i}\log F(t_1,t_2),~i=1,2$ are the components of bivariate reversed hazard rate;
\item[(ii)] the bivariate EIT is defined by the vector $m^X(t_1,t_2)=\left(m_1^X(t_1,t_2),m_2^X(t_1,t_2)\right)$ where $m_i^X(t_1,t_2)=E\left(t_i-X_i|X_1<t_1,X_2<t_2\right),~i=1,2$. For $i=1$,
    $$m_1^X(t_1,t_2)=\frac{1}{F(t_1,t_2)}\int_0^{t_1}F(x_1,t_2)dx_1,$$
    which measures the expected waiting time of the first component conditioned on the fact that both the components were failed before times $t_1$ and $t_2$, respectively.
\end{d1}
Note that, (\ref{eq3.1}) can alternatively be written as
\begin{equation}\label{eq3.2+}\overline\varepsilon^*_1(X;t_1,t_2)=m_1^X(t_1,t_2)\log F(t_1,t_2)-\int_0^{t_1}\frac{F(x_1,t_2)}{F(t_1,t_2)}\log F(x_1,t_2)dx_1,\end{equation}
where the above right hand side integral can be found to have a nice probabilistic meaning as follows. Let us set, for $0\leqslant a<b$,
$$T_1^{(2)}(a,b;t_2):=-\int_a^b \log F(x_1,t_2)dx_1.$$
Its partial derivative is closely related to the distribution function of $(X_1,X_2)$. Indeed, from above we have $\frac{\partial}{\partial t_1}T_1^{(2)}(a,t_1;t_2)=-\log F(t_1,t_2).$
Then,
\begin{eqnarray*}-\int_0^{t_1}\frac{F(x_1,t_2)}{F(t_1,t_2)}\log F(x_1,t_2)dx_1&=&-\frac{1}{F(t_1,t_2)}\int_0^{t_1}\left(\int_0^{x_1}f(u,t_2)du\right)\log F(x_1,t_2)dx_1\\
&=&-\frac{1}{F(t_1,t_2)}\int_0^{t_1}f(u,t_2)\left(\int_u^{t_1}\log F(x_1,t_2)dx_1\right)du\\
&=&E\left[T_1^{(2)}(X_1,t_1;t_2)|X_1<t_1,X_2<t_2\right]
\end{eqnarray*}
So, from (\ref{eq3.2+}) it can be written that
$$\overline\varepsilon^*_1(X;t_1,t_2)=m_1^X(t_1,t_2)\log F(t_1,t_2)+E\left[T_1^{(2)}(X_1,t_1;t_2)|X_1<t_1,X_2<t_2\right].$$
Similarly, (\ref{eq3.2}) can also be written as
\begin{eqnarray*}
\overline\varepsilon^*_2(X;t_1,t_2)&=&m_2^X(t_1,t_2)\log F(t_1,t_2)-\int_0^{t_2}\frac{F(t_1,x_2)}{F(t_1,t_2)}\log F(t_1,x_2)dx_2\\
&=&m_2^X(t_1,t_2)\log F(t_1,t_2)+E\left[T_2^{(2)}(X_2,t_2;t_1)|X_1<t_1,X_2<t_2\right],
\end{eqnarray*}
where $T_2^{(2)}(a,b;t_1)=-\int_a^b \log F(t_1,x_2)dx_2.$ Differentiating (\ref{eq3.1}) and (\ref{eq3.2}) with respect to $t_1$ and $t_2$, respectively we get in general
\begin{eqnarray}\label{eq3.3}
\frac{\partial}{\partial t_i}\overline\varepsilon^*_i(X;t_1,t_2)=\phi_i^X(t_1,t_2)\left[m_i^X(t_1,t_2)-\overline\varepsilon^*_i(X;t_1,t_2)\right],~~i=1,2.
\end{eqnarray}
Now we have the following theorem.
\begin{t1}\label{th3.1} For $t_1,t_2>0$, $\overline\varepsilon^*_i(X;t_1,t_2)$ is increasing in $t_i$, if and only if
$$\overline\varepsilon^*_i(X;t_1,t_2)\leqslant m_i^X(t_1,t_2),~~i=1,2.$$
\end{t1}
\begin{ex} Let $X$ follow the distribution as given in Example \ref{ex2.1}. Then, for $i=1,2$, $m_i^X(t_1,t_2)=t_i/2$ and $\overline\varepsilon^*_i(X;t_1,t_2)=t_i/4$. Here $\overline\varepsilon^*_i(X;t_1,t_2)$ is increasing in $t_i$ and satisfy the above inequality.
\end{ex}
\begin{r1}
Di Crescenzo and Longobardi (2009) pointed out that $\overline\varepsilon_{X}(t)$, defined in (\ref{eq1.5}), cannot be decreasing in $t$ for any random variable $X$ with support $(0,b)$ with $b$ finite or infinite. So, for $i=1$, if $t_2$ can be fixed at $b_2$, then $\overline\varepsilon^*_1(X;t_1,t_2)=\overline\varepsilon_{X_1}(t_1)$ cannot be decreasing in $t_1$. Thus for all $t_j>0$, it can be written that $\overline\varepsilon^*_i(X;t_1,t_2),~i=1,2,~i\neq j,$ cannot be decreasing in $t_i$.
\end{r1}
Although $\overline\varepsilon^*_i(X;t_1,t_2)$ is not decreasing in $t_i$ for all $t_j>0$, $i,j=1,2,~i\neq j$, the following example shows that $\overline\varepsilon^*_i(X;t_1,t_2)$ can be increasing in $t_i$ for all $t_j>0$.
\begin{ex}\label{ex3.1}
If $F(x_1,x_2)$ denotes joint distribution function of the random vector $X=(X_1,X_2)$, which follows bivariate extreme value distribution of type B, then for all $x_1,~x_2\geqslant 0$,
$$F(x_1,x_2)=exp\left[-\left(e^{-mx_1}+e^{-mx_2}\right)^{1/m}\right],~m\geqslant 1.$$
For properties of this distribution one may refer to Kotz et al. (2000). From Figure \ref{fig1} it is clear that $\overline\varepsilon^*_i(X;t_1,t_2)$ for this distribution taking $m=2$, is increasing in $t_i$ for all $t_j>0$, $i,j=1,2,~i\neq j$. It is to be mentioned here that while plotting curves, the substitutions $t_1=-\ln x$ and $t_2=-\ln y$ have been used so that
$\overline\varepsilon^*_1(X;-\ln x,-\ln y)=a_1(x,y)$ and $\overline\varepsilon^*_2(X;-\ln x,-\ln y)=a_2(x,y)$, say.
\begin{figure}[ht]
\centering
\begin{minipage}[b]{0.45\linewidth}
\includegraphics[height=7 cm]{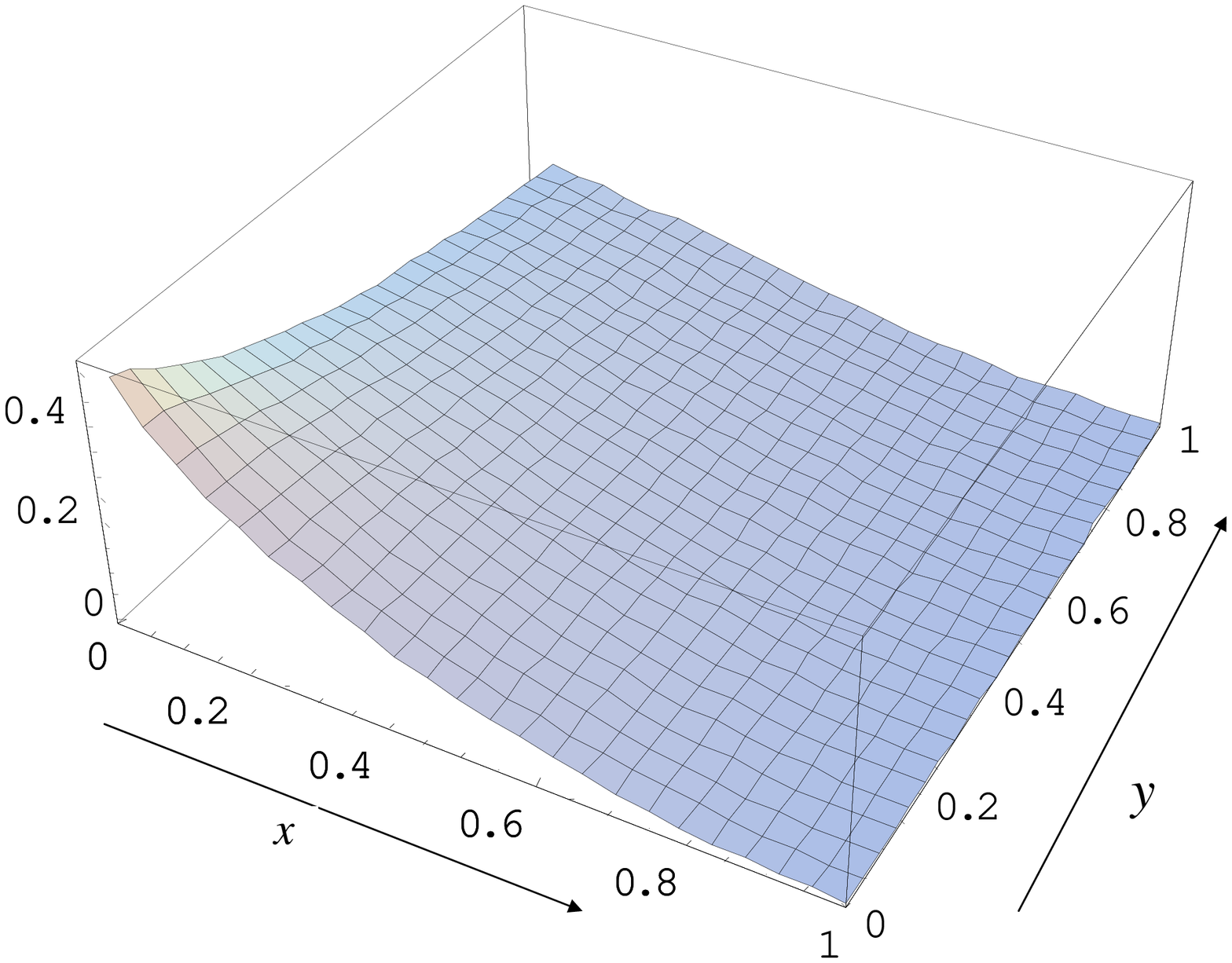}
\centering{Graph of $a_1(x,y)$}
\end{minipage}
\quad
\begin{minipage}[b]{0.45\linewidth}
\includegraphics[height=7 cm]{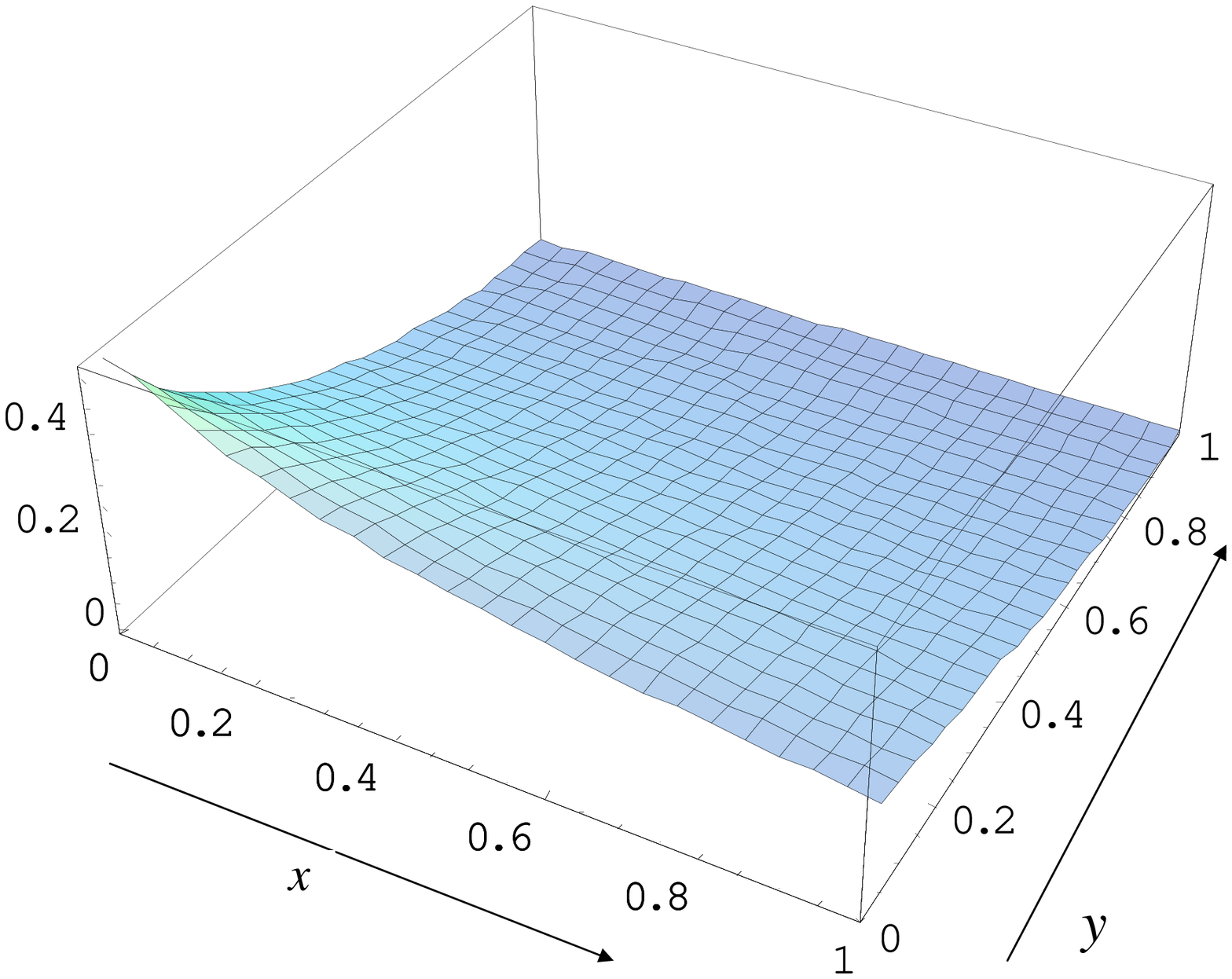}
\centering{Graph of $a_2(x,y)$}
\end{minipage}
\caption{\label{fig1}Graphical representations of $\left (a_1,a_2\right)$ (Example \ref{ex3.1})}
\end{figure}
\end{ex}
\hspace*{.2in} In the following theorem we obtain a functional relationship between the vector CDCPE and vector valued EIT. This relationship is useful in the sense that in many statistical models one may have information about bivariate EIT.
\begin{t1}\label{th3.2} Let $X=(X_1,X_2)$ be an absolutely continuous nonnegative bivariate random vector with finite $\overline\varepsilon^*_i(X;t_1,t_2)$ and bivariate EIT components $m_i^X(t_1,t_2),~i=1,2$. Then for all $t_1,t_2\geqslant 0$,
$$\overline\varepsilon^*_1(X;t_1,t_2)=\int_0^{t_1}m_1^X(x_1,t_2)f_1(x_1;t_1,t_2)dx_1$$
$${\rm and}~~\overline\varepsilon^*_2(X;t_1,t_2)=\int_0^{t_2}m_2^X(t_1,x_2)f_2(x_2;t_1,t_2)dx_2,$$
where $f_i(x_i;t_1,t_2)$ is the density function of $(X_i|X_1<t_1,X_2<t_2),~i=1,2.$
\end{t1}
Proof: Using by parts technique on the right hand side integral of (\ref{eq3.2+}), we get
$$\int_0^{t_1}\frac{F(x_1,t_2)}{F(t_1,t_2)}\log F(x_1,t_2)dx_1=m_1^X(t_1,t_2)\log F(t_1,t_2)-\int_0^{t_1}\frac{\frac{\partial F(x_1,t_2)}{\partial x_1}}{F(t_1,t_2)}\left(\int_0^{x_1}\frac{F(u,t_2)}{F(x_1,t_2)}du\right)dx_1.$$
Hence the first part follows from (\ref{eq3.2+}). The proof for second part is analogous. $\hfill\square$\\

\hspace*{.2in} On using Theorem \ref{th3.1} and \ref{th3.2}, we have the following result. The proof is omitted.
\begin{c1} For $t_1,t_2>0$, if $m_i^X(t_1,t_2)$ is increasing in $t_i,~i=1,2$ then $\overline\varepsilon^*_i(X;t_1,t_2)$ is also increasing in $t_i$.
\end{c1}

\hspace*{.2in} Let us set, for $0\leqslant a<b$,
$$\tau_1^{(2)}(a,b;t_2):=-\int_a^b \log \frac{F(x_1,t_2)}{F(t_1,t_2)}dx_1$$
$${\rm and}~~\tau_2^{(2)}(a,b;t_1):=-\int_a^b \log \frac{F(t_1,x_2)}{F(t_1,t_2)}dx_2.$$
Then an alternative expression to (\ref{eq3.1}) and (\ref{eq3.2}) is given hereafter. The proof is an immediate consequence of Fubini's theorem.
\begin{t1}\label{th3.3} Let $X=(X_1,X_2)$ be a nonnegative bivariate random vector with $\overline\varepsilon^*_i(X;t_1,t_2)$ finite. Then
$$\overline\varepsilon^*_1(X;t_1,t_2)=E\left[\tau_1^{(2)}(X_1,t_1;t_2)|X_1<t_1,X_2<t_2\right].$$
Similarly, $\overline\varepsilon^*_2(X;t_1,t_2)=E\left[\tau_2^{(2)}(X_2,t_2;t_1)|X_1<t_1,X_2<t_2\right]$.
\end{t1}
\hspace*{.2in} Let us discuss the effect of linear transformation on CDCPE. The proof is immediate from (\ref{eq3.1}) and (\ref{eq3.2}).
\begin{t1}\label{th3.4} Let $F$ and $G$ be the bivariate distribution functions of two nonnegative random vectors $X=(X_1,X_2)$ and $Y=(Y_1,Y_2)$, respectively where $Y_i=c_iX_i+d_i$ with $c_i>0$ and $d_i\geqslant0$ for $i=1,2$. Then
$$\overline\varepsilon^*_i(Y;t_1,t_2)=c_i\overline\varepsilon^*_i\left(X;\frac{t_1-d_1}{c_1},\frac{t_2-d_2}{c_2}\right),\qquad t_i\geqslant d_i,$$
which shows that CDCPE is a shift dependent measure.
\end{t1}
\begin{c1} Let $Y=(Y_1,Y_2)$ be a nonnegative bivariate random vector where $Y_i=c_iX_i+d_i$ with $c_i>0$ and $d_i\geqslant0$ for $i=1,2$. Then, $\overline\varepsilon^*_i(Y;t_1,t_2)$ is increasing in $t_i$ if and only if $\overline\varepsilon^*_i(X;t_1,t_2)$ is also so.
\end{c1}
\hspace*{.2in} The effect of monotonic transformation on CDCPE is given in the next theorem.
\begin{t1} Let $X=(X_1,X_2)$ be a nonnegative bivariate random vector. Suppose $\varphi(x)$ is strictly monotonic, continuous and differentiable function on $[0,\infty)$ with $\varphi(0)=0$ and derivative $\varphi'(x)$. If $\varphi$ is an increasing function with $a\leqslant\varphi'\leqslant b$, $a,b>0$ then
$$b\overline\varepsilon^*_i(X;\varphi^{-1}(t_1),\varphi^{-1}(t_2))\leqslant\overline\varepsilon^*_i\left(X_\varphi;t_1,t_2\right)
\leqslant a\overline\varepsilon^*_i(X;\varphi^{-1}(t_1),\varphi^{-1}(t_2)),~i=1,2,$$
where $X_\varphi=(\varphi(X_1),\varphi(X_2))$ is a bivariate random vector. If $\varphi$ is decreasing with $a\leqslant-\varphi'\leqslant b$, then
$$b\varepsilon_i(X;\varphi^{-1}(t_1),\varphi^{-1}(t_2))\leqslant\overline\varepsilon^*_i\left(X_\varphi;t_1,t_2\right)
\leqslant a\varepsilon_i(X;\varphi^{-1}(t_1),\varphi^{-1}(t_2)),$$
where $\varepsilon_i(X;t_1,t_2),~i=1,2$ is conditional DCRE as studied by Rajesh et al. (2014).
\end{t1}
Proof: If $\varphi$ is increasing then from (\ref{eq3.1}) and (\ref{eq3.2}), we have
$$\overline\varepsilon^*_1(X_\varphi;t_1,t_2)=-\int_0^{\varphi^{-1}(t_1)}\frac{F\left(u,\varphi^{-1}(t_2)\right)}{F\left(\varphi^{-1}(t_1),\varphi^{-1}(t_2)\right)}
\log\frac{F\left(u,\varphi^{-1}(t_2)\right)}{F\left(\varphi^{-1}(t_1),\varphi^{-1}(t_2)\right)}\varphi'(u)du$$
$${\rm and},~~\overline\varepsilon^*_2(X_\varphi;t_1,t_2)=-\int_0^{\varphi^{-1}(t_2)}\frac{F\left(\varphi^{-1}(t_1),v\right)}{F\left(\varphi^{-1}(t_1),\varphi^{-1}(t_2)\right)}
\log\frac{F\left(\varphi^{-1}(t_1),v\right)}{F\left(\varphi^{-1}(t_1),\varphi^{-1}(t_2)\right)}\varphi'(v)dv.$$
Hence the first part follows on using $a\leqslant\varphi'\leqslant b$. If $\varphi$ is decreasing we similarly obtain the second part of the proof. $\hfill\square$
\begin{t1} For a strictly increasing function $\varphi(\cdot)$ with $\varphi(0)=0$, $\overline\varepsilon^*_i(X_\varphi;t_1,t_2)$ is increasing in $t_i$ if $\overline\varepsilon^*_i(X;t_1,t_2)$ is increasing in $t_i$ for $i=1,2$ and $\varphi'<1$.
\end{t1}
Proof: Let $\varphi$ be strictly increasing. Then for $i=1$, we get after some algebraic calculation
\small\begin{eqnarray*}&&\frac{\partial}{\partial t_1}\left(\frac{\overline\varepsilon^*_1(X_\varphi;t_1,t_2)}{\overline\varepsilon^*_1(X;\varphi^{-1}(t_1),\varphi^{-1}(t_2))}\right)\\
&\stackrel{{\rm sign}}{=}&\int_0^{\varphi^{-1}(t_1)}F(u,\varphi^{-1}(t_2))du\int_0^{\varphi^{-1}(t_1)}F(u,\varphi^{-1}(t_2))\log F(u,\varphi^{-1}(t_2))\varphi'(u)du\\
&&-\int_0^{\varphi^{-1}(t_1)}F(u,\varphi^{-1}(t_2))\varphi'(u)du\int_0^{\varphi^{-1}(t_1)}F(u,\varphi^{-1}(t_2))\log F(u,\varphi^{-1}(t_2))du\\
&\geqslant&0,
\end{eqnarray*}
where the last inequality follows on using that $\varphi'<1$. Thus, $\overline\varepsilon^*_1(X_\varphi;t_1,t_2)$ is increasing in $t_1$ if $\overline\varepsilon^*_1(X;\varphi^{-1}(t_1),\varphi^{-1}(t_2))$, or equivalently, $\overline\varepsilon^*_1(X;t_1,t_2)$ is increasing in $t_1$. The result follows analogously for $i=2$. $\hfill\square$\\

\hspace*{.2in} Now we define the following stochastic orders between two bivariate random vectors based on CDCPE. For more on stochastic order one may refer to Shaked and Shanthikumar (2007).
\begin{d1} For two nonnegative bivariate random vectors $X=(X_1,X_2)$ and $Y=(Y_1,Y_2)$, $X$ is said to be greater than $Y$ in CDCPE order (written as $X\geqslant_{CDCPE}Y$) if $\overline\varepsilon^*_i(X;t_1,t_2)\leqslant\overline\varepsilon^*_i(Y;t_1,t_2)$ for all $t_i\geqslant0,~i=1,2$.
\end{d1}
\begin{r1}
It can be checked that the ordering defined above is reflexive, antisymmetric and transitive and thus partial ordering.
\end{r1}
\hspace*{.2in} Consider the following example to see that the ordering defined above is not equivalent to usual bivariate stochastic ordering.
\begin{ex}\label{ex3.2}
Let $X=\left(X_1,X_2\right)$ and $Y=\left(Y_1,Y_2\right)$ be two nonnegative bivariate random vector with distribution functions $F\left(t_1, t_2\right)$ and $G\left(t_1, t_2\right)$ respectively. Also let
$$F\left(t_1,t_2\right)=\frac{1}{\frac{4}{t_1}+\frac{4}{t_2}-1},~0\leqslant t_1\leqslant 4,~0\leqslant t_2\leqslant 4,$$
with marginals $F_{X_1}(t_1)=\frac{t_1}{4},~0\leqslant t_1\leqslant 4$ and $F_{X_2}(t_2)=\frac{t_2}{4},~0\leqslant t_2\leqslant 4$ and
$$G\left(t_1,t_2\right)=\frac{1}{\frac{1}{t_1^2}+\frac{1}{t_2^2}-1},~0\leqslant t_1\leqslant 1,~0\leqslant t_2\leqslant 1,$$
with marginals $G_{Y_1}(t_1)=t_1^2,~0\leqslant t_1\leqslant 1$ and $G_{Y_2}(t_2)=t_2^2,~0\leqslant t_2\leqslant 1$. From Figure \ref{fig2} it is clear that for all $0\leqslant t_1\leqslant 1$ and $0\leqslant t_2\leqslant 1$ and for $i=1,2$, $\overline\varepsilon^*_i(X;t_1,t_2)\geqslant\overline\varepsilon^*_i(Y;t_1,t_2)$, giving $X\leqslant_{CDCPE}Y$. Again $F_{X_1}(0.4)-G_{Y_1}(0.4)=-0.44$ and $F_{X_1}(0.1)-G_{Y_1}(0.1)=0.015$ proving $X_1\ngeqslant_{st} Y_1$. So, by Theorem 6.B.16 (c) of Shaked and Shanthikumar (2007) it can be concluded that $X\ngeqslant_{st} Y$.
\begin{figure}[ht]
\centering
\begin{minipage}[b]{0.45\linewidth}
\includegraphics[height=7 cm]{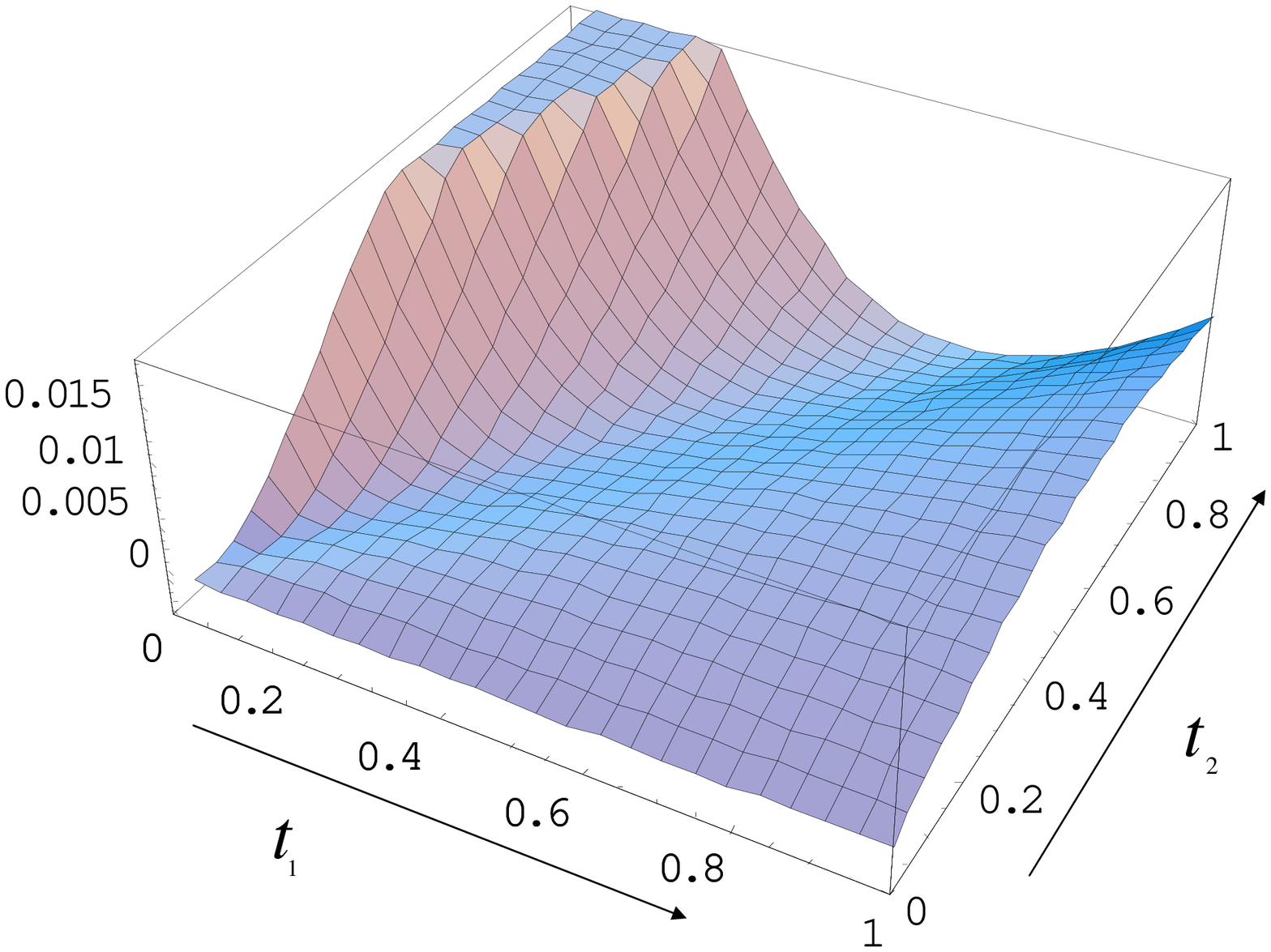}
\centering{Graph of $\overline\varepsilon^*_1(X;t_1,t_2)-\overline\varepsilon^*_1(Y;t_1,t_2)$}
\end{minipage}
\quad
\begin{minipage}[b]{0.45\linewidth}
\includegraphics[height=7 cm]{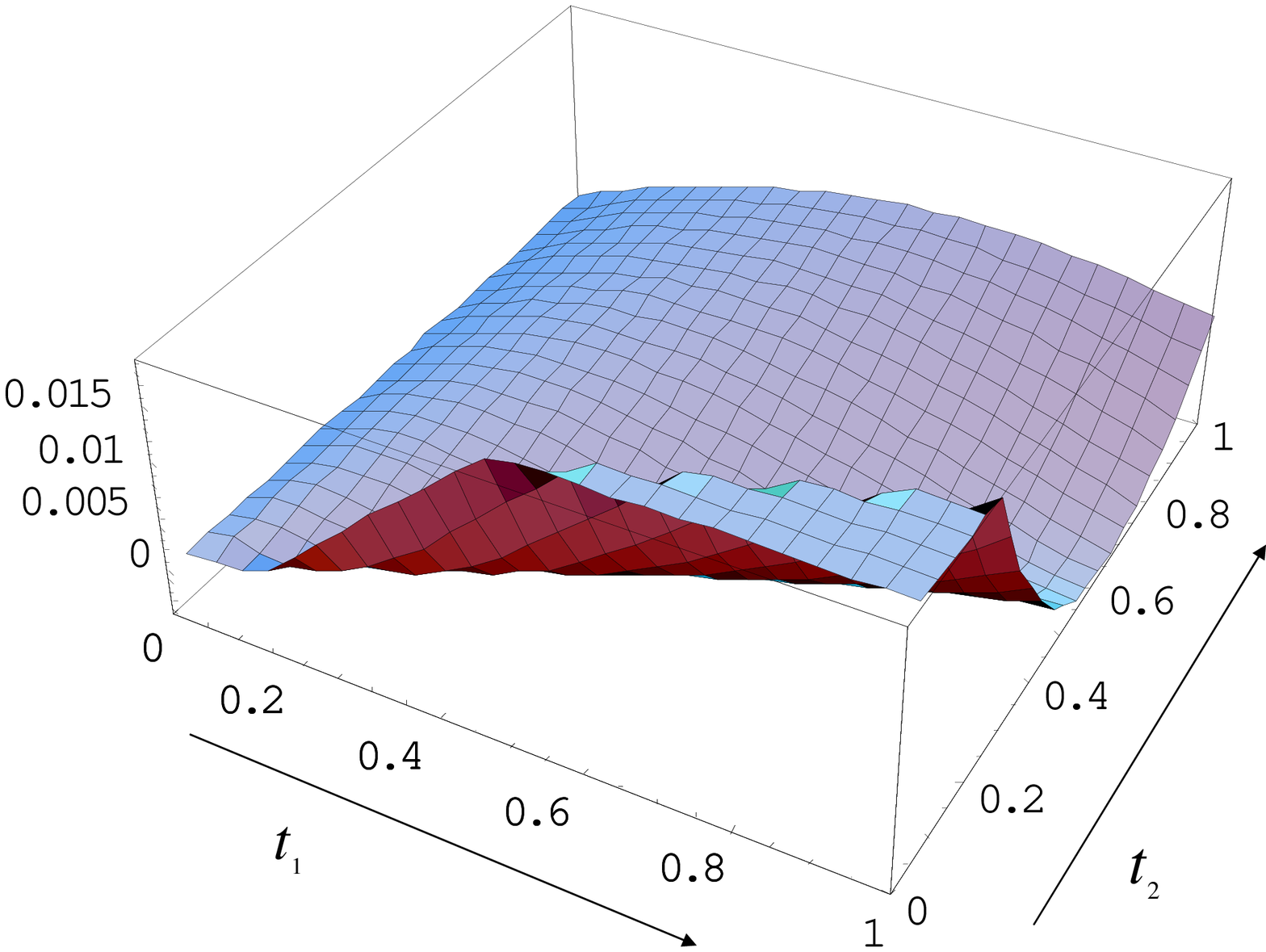}
\centering{Graph of $\overline\varepsilon^*_2(X;t_1,t_2)-\overline\varepsilon^*_2(Y;t_1,t_2)$}
\end{minipage}
\caption{\label{fig2}Graphical representations of $\overline\varepsilon^*_i(X;t_1,t_2)-\overline\varepsilon^*_i(Y;t_1,t_2)$ (Example \ref{ex3.2})}
\end{figure}
\end{ex}
\hspace*{.2in} The following result can be obtained from Theorem \ref{th3.4}.
\begin{t1}\label{th3.7} For two nonnegative bivariate random vectors $X=(X_1,X_2)$ and $X'=(X_1',X_2')$, let $Y_i=c_iX_i+d_i$ and $Y_i'=c_iX_i'+d_i$ with $c_i>0$ and $d_i\geqslant0$ for $i=1,2$. Then $Y\geqslant_{CDCPE}Y'$ if $X\geqslant_{CDCPE}X'$ where $Y'=(Y_1',Y_2')$.
\end{t1}
\hspace*{.2in} Now we have the more general result.
\begin{t1}\label{th3.8} Let $F,F',G$ and $G'$ be the joint distribution functions of bivariate random vectors $X,X',Y$ and $Y'$, respectively. Also let for $i=1,2$, $Y_i=a_iX_i+c_i$ and $Y_i'=a_iX'_i+d_i$ with $a_i>0$ and $d_i\geqslant c_i>0$. Then, $Y\leqslant_{CDCPE}Y'$ provided $X\leqslant_{CDCPE}X'$ and either $\overline\varepsilon^*_i(X;t_1,t_2)$ or $\overline\varepsilon^*_i(X';t_1,t_2)$ is increasing in $t_1$ as well as $t_2$.
\end{t1}
Proof: If $X\leqslant_{CDCPE}X'$, then using Theorem \ref{th3.4} it can be written that
\begin{eqnarray*}\overline\varepsilon^*_i(Y;t_1,t_2)&\geqslant& a_i\overline\varepsilon^*_i\left(X';\frac{t_1-c_1}{a_1},\frac{t_2-c_2}{a_2}\right)\\
&\geqslant& a_i\overline\varepsilon^*_i\left(X';\frac{t_1-d_1}{a_1},\frac{t_2-d_2}{a_2}\right)=\overline\varepsilon^*_i(Y';t_1,t_2),
\end{eqnarray*}
where the last inequality follows on using the fact that $\overline\varepsilon^*_i(X';t_1,t_2)$ is increasing in $t_1$ and $t_2$. Similarly, the result follows if $\overline\varepsilon^*_i(X;t_1,t_2)$ is increasing in $t_1$ and $t_2$. $\hfill\square$

\hspace*{.2in} In the recent past, the researchers have shown more interest in characterization of distributions. An important question regarding the CDCPE is whether it characterizes the underlying distribution function uniquely. In the following theorem we show that $\overline\varepsilon^*_i(X;t_1,t_2),~i=1,2$ determines the distribution function uniquely.
\begin{t1}
Let $X=\left(X_1, X_2\right)$ be a nonnegative bivariate random variable having absolutely continuous distribution function $F$ with respect to the Lebesgue measure. Then CDCPE of $X$, defined in (\ref{eq3.1}) and (\ref{eq3.2}), uniquely determines the distribution function provided they are finite.
\end{t1}
Proof: Let $X$ and $Y$ be two bivariate random variables having joint distribution functions $F$ and $G$, respectively. Also let, for all $t_1,t_2\geqslant 0$,
$$\overline\varepsilon^*_i(X;t_1,t_2)=\overline\varepsilon^*_i(Y;t_1,t_2),\ i=1,\ 2.$$
Differentiating $\overline\varepsilon^*_i(X;t_1,t_2)$ and $\overline\varepsilon^*_i(Y;t_1,t_2)$ with respect to $t_i,~i=1,2$ and, on using the relation $\phi_i^X(t_1,t_2)m_i^X(t_1,t_2)=1-\frac{\partial}{\partial t_i}m_i^X(t_1,t_2),$ we have from (\ref{eq3.3})
$$\frac{\partial}{\partial t_i}\left(m_i^X(t_1,t_2)\right)=\frac{\frac{\partial}{\partial t_i}\overline\varepsilon^*_i(X;t_1,t_2)m_i^X(t_1,t_2)+\overline\varepsilon^*_i(X;t_1,t_2)-m_i^X(t_1,t_2)}{\overline\varepsilon^*_i(X;t_1,t_2)-m_i^X(t_1,t_2)}$$
and
$$\frac{\partial}{\partial t_i}\left(m_i^Y(t_1,t_2)\right)=\frac{\frac{\partial}{\partial t_i}\overline\varepsilon^*_i(Y;t_1,t_2)m_i^Y(t_1,t_2)+\overline\varepsilon^*_i(Y;t_1,t_2)-m_i^Y(t_1,t_2)}{\overline\varepsilon^*_i(Y;t_1,t_2)-m_i^Y(t_1,t_2)}$$
respectively. Suppose that
$$\overline\varepsilon^*_i(X;t_1,t_2)=\overline\varepsilon^*_i(Y;t_1,t_2)=\theta_i\left(\bf t\right)$$
and
$$\psi_i\left({\bf t},z\right)=\frac{\frac{\partial}{\partial t_i}\theta_i\left({\bf t}\right)z+\theta_i\left({\bf t}\right)-z}{\theta_i\left({\bf t}\right)-z},~~{\bf t}=\left(t_1,t_2\right),$$
for $i=1,2$. Thus, we can write
$$\frac{\partial}{\partial t_i}\left(m_i^X(t_1,t_2)\right)=\psi_i\left({\bf t},m_i^X(t_1,t_2)\right)$$
and
$$\frac{\partial}{\partial t_i}\left(m_i^Y(t_1,t_2)\right)=\psi_i\left({\bf t},m_i^Y(t_1,t_2)\right).$$
In a recent paper, Thapliyal et al. (2013) proved that dynamic cumulative past entropy determines the distribution function uniquely, which in turn uniquely determines expected inactivity time. Thus $\overline\varepsilon^*_i(X;t_1,t_2)$ determines $m_i^X(t_1,t_2)$, $i=1,2$. Again, using the fact that vector valued expected inactivity time uniquely determines the bivariate distribution function, the proof is complete. $\hfill\square$\\

\hspace*{.2in} Now we study the characterization result for uniform distribution. The proof follows easily as $X_1$ and $X_2$ are independent.
\begin{t1}
Let $X=(X_1,X_2)$ be a bivariate random variable with distribution function $F$. Then $X$ is said to follow bivariate uniform distribution with distribution function
$$F\left(t_1,t_2\right)=\frac{t_1t_2}{b~d},~0\leqslant t_1\leqslant b,~ 0\leqslant t_2\leqslant d,$$
if and only if $\overline\varepsilon^*_i(X;t_1,t_2)=\frac{t_i}{4},~i=1,2$.
\end{t1}

\hspace*{.2in} In the following theorem we give another characterization result with dependent components.
\begin{t1} Let $X$ be a nonnegative bivariate random vector with $\overline\varepsilon^*_i(X;t_1,t_2)$ finite and the components of bivariate EIT are $m_i^X(t_1,t_2),~i=1,2$ for all $t_i\geqslant0$. Then, for $0<t_1,t_2<1,~\theta\leqslant0$,
\begin{equation}\label{eq3.4}
\overline\varepsilon^*_i(X;t_1,t_2)=\left(\frac{1+\theta\log t_j}{2+\theta\log t_j}\right)m_i^X(t_1,t_2),~i,j=1,2,~i\neq j,
\end{equation}
if and only if $X$ is distributed as bivariate uniform with
\begin{equation}\label{eq3.5}
F(t_1,t_2)=t_1^{1+\theta\log t_2}t_2,~0<t_1,t_2<1,~\theta\leqslant0.
\end{equation}
\end{t1}
Proof: The if part is trivial on noting that if $X$ follows the distribution (\ref{eq3.5}) then
$$m_i^X(t_1,t_2)=\frac{t_i}{2+\theta\log t_j},~~~{\rm and}~~~\overline\varepsilon^*_i(X;t_1,t_2)=\frac{t_i(1+\theta\log t_j)}{(2+\theta\log t_j)^2},~~i,j=1,2,~i\neq j.$$
To prove the converse let us assume that (\ref{eq3.4}) holds. Then, differentiating (\ref{eq3.4}) with respect to $t_i$, and using (\ref{eq3.3}), we get after some algebraic manipulation
$$\frac{\partial}{\partial t_i}m_i^X(t_1,t_2)=\frac{1}{2+\theta\log t_j},~i,j=1,2,~i\neq j$$
which on integration gives
$$m_i^X(t_1,t_2)=\frac{t_i}{2+\theta\log t_j}+c_i(t_j),$$
where $c_i(t_j)$ is a constant of integration. Now, $c_i(t_j)=0$ as $m_i^X(t_1,t_2)\rightarrow0$ for $t_i\rightarrow0$, which in turn gives the bivariate EIT of (\ref{eq3.5}). Hence the result follows on using the fact that bivariate EIT determines the distribution uniquely. $\hfill\square$\\

\hspace*{.2in} The following theorem gives a characterization of the bivariate power distribution. This result extends Theorem 6.2 of Di Crescenzo and Longobardi (2009) to bivariate setup.
\begin{t1} Let $X$ be a nonnegative bivariate random vector in the support $(0,b_1)\times(0,b_2)$, $b_i<\infty,~i=1,2$ with $\overline\varepsilon^*_i(X;t_1,t_2)$ finite. Then
\begin{equation}\label{eq3.6}
\overline\varepsilon^*_i(X;t_1,t_2)=c_i(t_j)m_i^X(t_1,t_2),~i,j=1,2,~i\neq j,
\end{equation}
where $c_i(t_j)\in(0,1)$ is a function independent of $t_i$, characterizes the bivariate power distribution
\begin{equation}\label{eq3.7}
F(t_1,t_2)=\left(\frac{t_1}{b_1}\right)^{c_1}\left(\frac{t_2}{b_2}\right)^{c_2+\theta\log\left(\frac{t_1}{b_1}\right)},~\theta\leqslant0,
\end{equation}
where $c_i=c_i(b_j)/\left[1-c_i(b_j)\right]$.
\end{t1}
Proof: The if part is straightforward. to prove the reverse implication, if (\ref{eq3.6}) holds, then differentiating both the sides with respect to $t_i$, and on using (\ref{eq3.3}), we get after some algebraic manipulation
$$\frac{\partial}{\partial t_i}m_i^X(t_1,t_2)=\left[1-c_i(t_j)\right],~i,j=1,2,~i\neq j$$
which on integration gives
$$m_i^X(t_1,t_2)=\left[1-c_i(t_j)\right]t_i+k_i(t_j),$$
where $k_i(t_j)$ is a constant of integration. Now, $k_i(t_j)=0$ as $m_i^X(t_1,t_2)\rightarrow0$ for $t_i\rightarrow0$. The rest of the proof follows from Theorem 2.1 of Nair and Asha (2008).

\subsection{Conditional dynamic CPE for $X_i$ given $X_j=t_j$}
\begin{figure}[ht]
\centering
\begin{minipage}[b]{0.45\linewidth}
\includegraphics[height=7 cm]{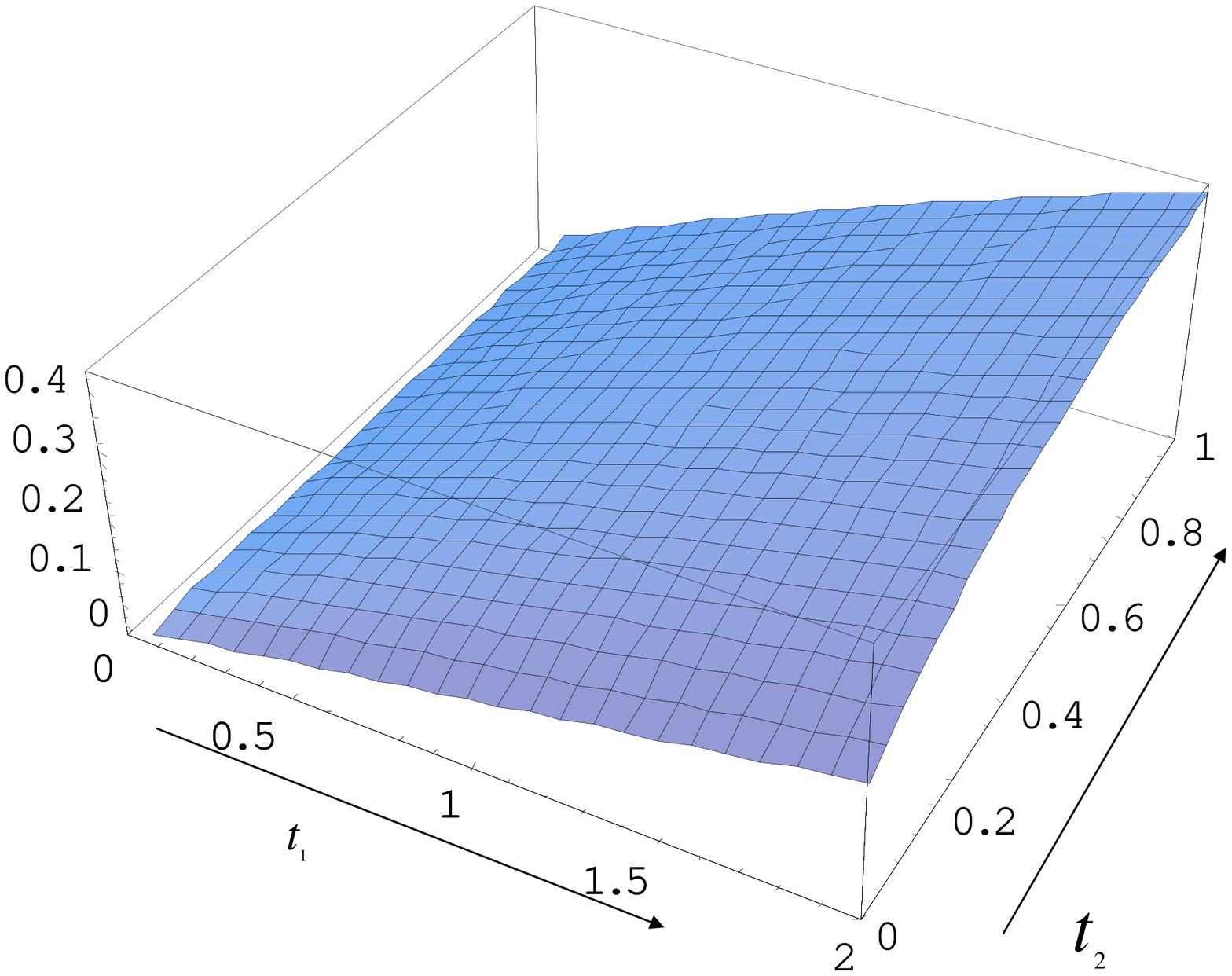}
\centering{Graph of $\overline m_1^X\left(t_1|t_2\right)-\overline\gamma^*_1(X;t_1,t_2)$}
\end{minipage}
\quad
\begin{minipage}[b]{0.45\linewidth}
\includegraphics[height=7 cm]{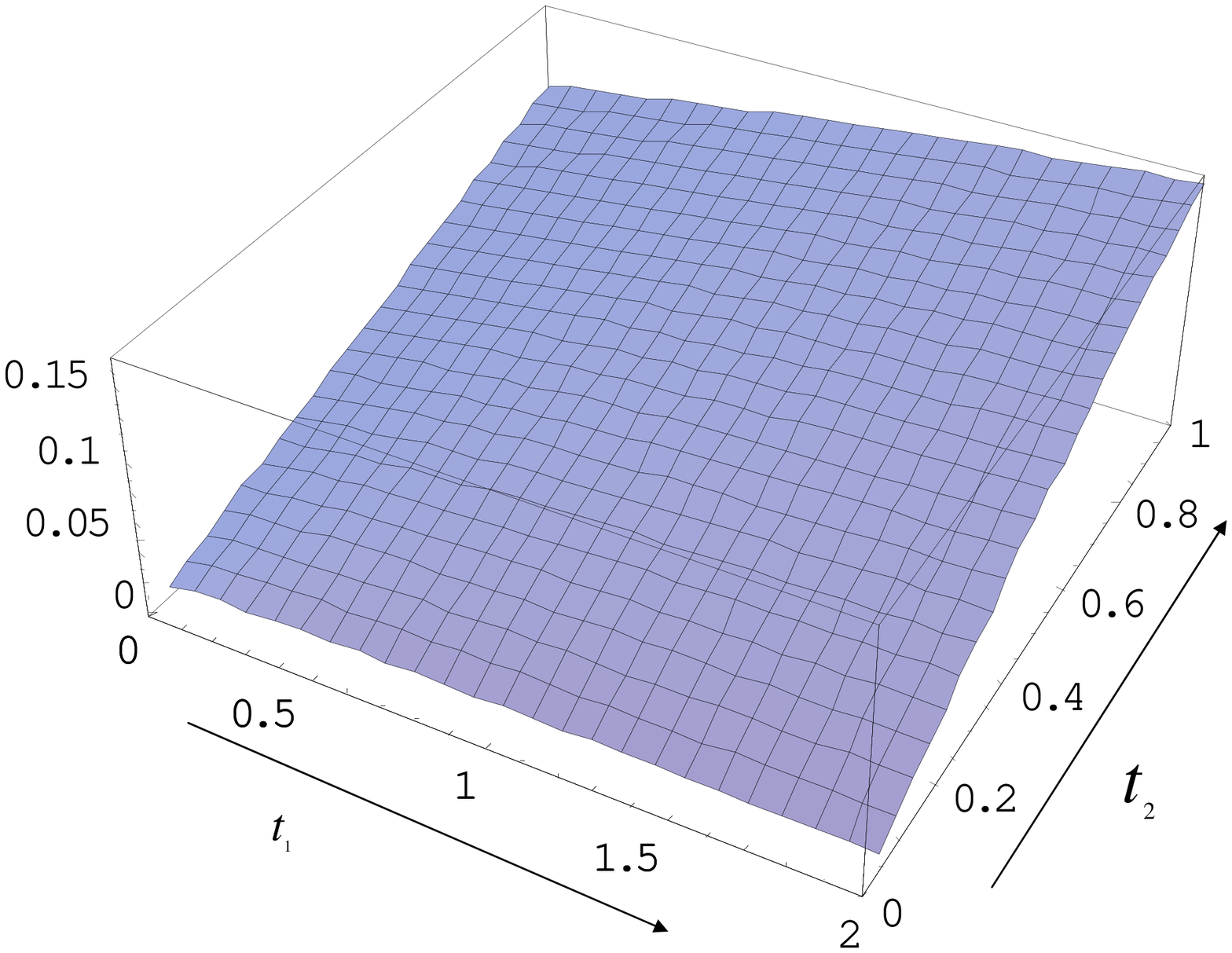}
\centering{Graph of $\overline m_2^X\left(t_2|t_1\right)-\overline\gamma^*_2(X;t_1,t_2)$}
\end{minipage}
\caption{\label{fig3}Graph of $\overline m_i^X\left(t_i|t_j\right)-\overline\gamma^*_i(X;t_1,t_2)$ (Example \ref{ex3.3})}
\end{figure}

The determination of the joint distribution function of $X=(X_1,X_2)$, when conditional distributions of $(X_1|X_2=t_2)$ and $(X_2|X_1=t_1)$ are known, has been an important problem dealt with by many researchers in the past. This approach of identifying a bivariate density using the conditionals is called the conditional specification of the joint distribution (see Arnold et al., 1999). These conditional models are often useful in many two component reliability systems, when the operational status of one component is known. Let the distribution function of $\widetilde{Y_i}^*=(X_i|X_i<t_i,X_j=t_j),$ $i,j=1,2$, $i\neq j$ be $F^*_i\left(t_i|t_j\right)$. Then, for an absolutely continuous nonnegative bivariate random vector $X$, the conditional dynamic CPE of $\widetilde{Y_i}^*$ is defined as
\begin{equation}\label{eq4.1}
\overline\gamma^*_i(X;t_1,t_2)=-\int_0^{t_i}\frac{F^*_i(x_i|t_j)}{F^*_i(t_i|t_j)}\log\frac{F^*_i(x_i|t_j)}{F^*_i(t_i|t_j)}dx_i,\qquad x_i<t_i,
\end{equation}
$i,j=1,2$, $i\neq j$. In particular, if $X_1$ and $X_2$ are independent, then (\ref{eq4.1}) reduces to marginal dynamic CPE of $X_i,~i=1,2$ as given in (\ref{eq1.5}). Following Roy (2002) the bivariate reversed hazard rate of $X=(X_1,X_2)$ is also defined by a vector, $\overline\phi^X\left(t_i|t_j\right)=\left(\overline\phi_1^X(t_1|t_2),\overline\phi_2^X(t_2|t_1)\right)$, where $\overline\phi_i^X\left(t_i|t_j\right)=\frac{\partial}{\partial t_i}\log F^*_i\left(t_i|t_j\right)$, $i,j=1,2$, $i\neq j$. For $i=1$, $\overline\phi_1^X(t_1|t_2)\Delta t_1$ is the probability of failure of the first component in the interval $(t_1-\Delta t_1, t_1]$ given that it has failed before $t_1$ and the failure time of the second is $t_2$.
Another definition of bivariate EIT of $X=(X_1,X_2)$ is given by Kayid (2006) as a vector, $\overline m^X\left(t_i|t_j\right)=\left(\overline m_1^X(t_1|t_2),\overline m_2^X(t_2|t_1)\right)$, where $\overline m_i^X\left(t_i|t_j\right)=E\left(t_i-X_i|X_i<t_i, X_j=t_j\right)$, $i,j=1,2$, $i\neq j$. For $i=1$,
$$\overline m_1^X(t_1|t_2)=\frac{1}{F^*_1(t_1|t_2)}\int_0^{t_1}F^*_1\left(x_1|t_2\right)dx_1,$$
which measures the expected waiting time of $X_1$ given that $X_1<t_1$ and $X_2=t_2$. Unlike $\phi^X(t_1,t_2)$ and $m^X(t_1,t_2)$, $\overline m^X\left(t_i|t_j\right)$ determines the distribution uniquely. But, $\overline\phi^X\left(t_i|t_j\right)$ does not provide $F(t_1,t_2)$ uniquely.\\
\hspace*{.2in} Differentiating (\ref{eq4.1}) with respect to $t_i$ and simplifying, we get
$$\frac{\partial}{\partial t_i}\overline\gamma^*_i(X;t_1,t_2)=\overline\phi_i^X\left(t_i|t_j\right)\left[\overline m_i^X\left(t_i|t_j\right)-\overline\gamma^*_i(X;t_1,t_2)\right],\ i,j=1,2, i\neq j.$$
Now we have the following theorem.
\begin{t1}\label{th3.13} For $t_1,t_2>0$, $\overline\gamma^*_i(X;t_1,t_2)$ is increasing in $t_i$, if and only if $$\overline\gamma^*_i(X;t_1,t_2)\leqslant\overline m_i^X\left(t_i|t_j\right), \ i,j=1,2, i\neq j.$$
\end{t1}
\hspace*{.2in} The following example gives an application of the above theorem.
\begin{ex}\label{ex3.3}
If the density function of a continuous bivariate random vector $X=\left(X_1,X_2\right)$, $f\left(x_1,x_2\right)$ is given by
$$f(x_1,x_2)= \left\{\begin{array}{ll}
    \frac{1}{6}(x_1+4x_2), & 0\leqslant x_1\leqslant 2,~ 0\leqslant x_2\leqslant 1 \\
                   0,                      & elsewhere,
                          \end{array}\right. $$
then it can be checked that for $i,j=1,2$ and $i\neq j$, $\overline\gamma^*_i(X;t_1,t_2)$ is increasing in $t_i$ for all $t_j\geqslant 0$. Again, Figure \ref{fig3} shows that for $0\leqslant t_1\leqslant 2$ and $0\leqslant t_2\leqslant 1$ and $i=1,2$, $i\neq j$, $\overline m_i^X\left(t_i|t_j\right)-\overline\gamma^*_i(X;t_1,t_2)$ are always positive, satisfying Theorem \ref{th3.13}.
\end{ex}
\hspace*{.2in} The next theorem, analogous to Theorem \ref{th3.2}, establishes relation between dynamic CPE and EIT of $X_i$ given $X_j=t_j$, $i,j=1,2$, $i\neq j$. The proof is omitted.
\begin{t1} For $t_1, t_2>0$, let $\overline\gamma^*_i(X;t_1,t_2)$ and $\overline m_i^X\left(t_i|t_j\right)$ are CDCPE and the components of bivariate EIT as defined above. If $f^*_i\left(t_i|t_j\right)$ denotes the density function of $\left(X_i|X_i<t_i,X_j=t_j\right)$, then
$$\overline\gamma^*_i(X;t_1,t_2)=\int_0^{t_i}\overline m_i^X\left(x_i|t_j\right)f^*_i\left(x_i|t_j\right)dx_i,\ i,j=1,2,\ i\neq j.$$
\end{t1}
\hspace*{.2in} Let $f\left(t_1,t_2\right)$ be the joint pdf of the bivariate random variable $X$, then (\ref{eq4.1}) can alternatively be written as
\begin{equation*}
\overline\gamma^*_i(X;t_1,t_2)=\overline m_i^X\left(t_i|t_j\right)\log\int_0^{t_i} f\left(z,t_j\right)dz-\int_0^{t_i}\frac{F^*_i\left(s|t_j\right)}{F^*_i\left(t_i|t_j\right)}\log\int_0^{s} f\left(z,t_j\right)dz\ ds,
\end{equation*}
$i,j=1,2,\ i\neq j$. So, taking
$$\overline T_i^{(2)}(a,b)=-\int_a^b \log\int_0^s f\left(z,t_j\right)dz\ ds,\ i,j=1,2,\ i\neq j,$$
and proceeding in the similar way as of previous subsection, it can be shown that
$$\overline\gamma^*_i(X;t_1,t_2)=\overline m_i^X\left(t_i|t_j\right)\log\int_0^{t_i} f\left(z,t_j\right)dz+E\left[\overline T_i^{(2)}(X_i,t_i)|X_i<t_i,X_j=t_j\right], \ i,j=1,2,\ i\neq j.$$
Again, for $0\leqslant a\leqslant b$, defining
$$\overline \tau_i^{(2)}(a,b)=-\int_a^b \log\frac{\int_0^s f\left(z,t_j\right)dz}{\int_0^{t_i} f\left(z,t_j\right)dz}\ ds,\ i,j=1,2,\ i\neq j,$$
as before, we have the following theorem, which is analogous to Theorem \ref{th3.3}.
\begin{t1} Let $X=\left(X_1,X_2\right)$ be a nonnegative bivariate random vector with $\overline\gamma^*_i(X;t_1,t_2)<\infty$, $i=1,2$. Then
$$\overline\gamma^*_i(X;t_1,t_2)=E\left[\overline \tau_i^{(2)}(X_i,t_i)|X_i<t_i,X_j=t_j\right],\ i,j=1,2,\ i\neq j.$$
\end{t1}
\hspace*{.2in} The effect of linear transformation on $\overline\gamma^*_i(X;t_1,t_2)$ is given in the following theorem.
\begin{t1}\label{th3.15} Let $X=\left(X_1,X_2\right)$ and $Y=\left(Y_1,Y_2\right)$ be two nonnegative bivariate random variables having distribution functions $F$ and $G$ respectively, where $Y_i=c_iX_i+d_i$, $i=1,2$. Then
$$\overline\gamma^*_i(Y;t_1,t_2)=c_i\overline\gamma^*_i\left(X;\frac{t_1-d_1}{c_1},\frac{t_2-d_2}{c_2}\right),\ t_i\geqslant d_i.$$
\end{t1}
\begin{c1}
Let $Y=(Y_1,Y_2)$ be a nonnegative bivariate random vector where $Y_i=c_iX_i+d_i$ with $c_i>0$ and $d_i\geqslant0$ for $i=1,2$. Then, $\overline\gamma^*_i(Y;t_1,t_2)$ is increasing in $t_i$ if and only if $\overline\gamma^*_i(X;t_1,t_2)$ is also so.
\end{c1}
The immediate consequence of Theorem \ref{th3.15} are the following two theorems. The proof being analogous to Theorems \ref{th3.7} and \ref{th3.8} is omitted.
\begin{t1} Let $X=(X_1,X_2)$ and $X'=(X_1',X_2')$ be two nonnegative bivariate random variables. Also let $Y=(Y_1,Y_2)$ and $Y'=(Y_1',Y_2')$ be two nonnegative bivariate random vectors such that $Y_i=c_i X_i+d_i$ and $Y'_i=c_iX'_i+d_i$ with $c_i>0$, $d_i\geqslant 0$ for $i=1,2$. If $\overline\gamma^*_i(X;t_1,t_2)\geqslant \overline\gamma^*_i(X';t_1,t_2)$ then $\overline\gamma^*_i(Y;t_1,t_2)\geqslant \overline\gamma^*_i(Y';t_1,t_2)$, for all $t_1, t_2\geqslant 0$, where $F,F',G,G'$ are the distribution functions of $X,X',Y$ and $Y'$, respectively.
\end{t1}
\begin{t1} Let $F,F',G,G'$ are the distribution functions of nonnegative bivariate random variables $X,X',Y$ and $Y'$, respectively. Also let  for $i=1,2$, $Y_i=a_i X_i+c_i$ and $Y'_i=a_iX'_i+d_i$ with $a_i>0$ and $d_i\geqslant c_i>0$. Then, for all $t_1, t_2\geqslant 0$, $\overline\gamma^*_i(Y;t_1,t_2)\geqslant \overline\gamma^*_i(Y';t_1,t_2)$ provided $\overline\gamma^*_i(X;t_1,t_2)\geqslant \overline\gamma^*_i(X';t_1,t_2)$ and either $\overline\gamma^*_i(X;t_1,t_2)$ or $\overline\gamma^*_i(X';t_1,t_2)$ is increasing in $t_1$ as well as $t_2$.
\end{t1}
\section*{Acknowledgements}
The authors thank the anonymous reviewers for their constructive comments which have considerably improved the content of the paper. We also thank Prof. Asok K. Nanda, IISER Kolkata for helpful discussions while revising the manuscript. The financial support (vide no. PSW-103/13-14 (ERO), ID No. WK4-031 dated 18.03.2014) from the University Grants Commission, Government of India, is acknowledged with thanks by Amarjit Kundu. Also Chanchal Kundu acknowledges with thanks the financial support (Ref. No. SR/FTP/MS-016/2012) rendered by the Department of Science and Technology, Government of India for doing this research.


\begin{thebibliography}{99}
{\small
\bibitem{a} Abbasnejad, M. (2011), Some Characterization results based on
dynamic survival and failure entropies. {\it Communications of the Korean Statistical Society},
{\bf18}, 1-12.
\bibitem{a1} Abbasnejad, M., Arghami, N.R, Morgenthaler, S. and Borzadaran, G.R.M (2010), On the
dynamic survival entropy. {\it Statistics and Probability Letters}, {\bf80}, 1962-1971.
\bibitem{acs1} Arnold, B.C., Castillo, E. and Sarabia, J.M. (1999), {\it Conditional Specification of Statistical Models}. Springer Verlag, New York.
\bibitem{acs2} Arnold, B.C., Castillo, E. and Sarabia, J.M. (2001), Conditionally specified distributions: an introduction (with discussion). {\it Statistical Science}, {\bf16(3)}, 249-274.
\bibitem{az} Asadi, M. and Zohrevand, Y. (2007), On the dynamic
cumulative residual entropy. {\it Journal of Statistical Planning
and Inference}, {\bf 137}, 1931-1941.
\bibitem{cb} Chamany, A. and Baratpour, S. (2014), A dynamic discrimination information based on cumulative residual entropy and its properties.
{\it Communications in Statistics- Theory \& Methods}, {\bf43(6)}, 1041-1049.
\bibitem{dl1} Di Crescenzo, A. and Longobardi, M. (2009), On cumulative entropies.
{\it Journal of Statistical Planning and Inference}, {\bf 139}, 4072-4087.
\bibitem{dl2} Di Crescenzo, A. and Longobardi, M. (2013), Stochastic comparisons of cumulative entropies.
{\it Stochastic Orders in Reliability and Risk}, Lecture Notes in Statistics 208, 167-182. Springer, New York.
\bibitem{hn} Hitha, N. and Nair, N.U. (1991), Characterizations of bivariate Lomax and finite range distribution. {\it Calcutta Statistical Association Bulletin}, {\bf41}, 163-167.
    \bibitem{k} Kayid, M. (2006), Multivariate mean inactivity time functions with reliability applications. {\it International Journal of Reliability and Applications}, {\bf7(2)}, 125-139.
\bibitem{ko} Kotz, S., Balakrishnan, N. and Johnson, N.L. (2000), {\it Continuous Multivariate Distributions: Models and Applications}. John Wiley \& Sons, Inc.
\bibitem{krm} Khorashadizadeh, M., Rezaei Roknabadi, A.H. and Mohtashami Borzadaran, G.R. (2013),
Doubly truncated (interval) cumulative residual and past entropy.
{\it Statistics and Probability Letters}, {\bf83}, 1464-1471.
\bibitem{kt} Kumar, V. and Taneja, H.C. (2011), Some characterization
results on generalized cumulative residual entropy measure.
{\it Statistics and Probability Letters}, {\bf81}, 1072-1077.
\bibitem{na} Nair, N.U. and Asha, G. (2008), Some characterizations based on bivariate reversed mean residual life. {\it ProbStat Forum}, {\bf1}, 1-14.
\bibitem{nn} Nair, N.U. and Nair, V.K.R (1988), A characterization of the bivariate exponential distribution. {\it Biometrical Journal}, {\bf30(1)}, 107-112.
\bibitem{ns} Navarro, J. and Sarabia, J.M. (2013), Reliability properties of bivariate conditional proportional hazard rate models. {\it Journal of Multivariate Analysis}, {\bf113}, 116-127.
\bibitem{naa} Navarro, J., del Aguila, Y. and Asadi, M. (2010),
Some new results on the cumulative residual entropy. {\it Journal of Statistical
Planning and Inference}, {\bf 140}, 310-322.
\bibitem{nsl1} Navarro, J., Sunoj, S.M and Linu, M.N. (2011), Characterizations of bivariate models using dynamic Kullback-Leibler discrimination measures. {\it Statistics and Probability Letters}, {\bf81 (11)}, 1594-1598.
\bibitem{nsl2} Navarro, J., Sunoj, S.M and Linu, M.N. (2014), Characterizations of bivariate models using some dynamic conditional information divergence measures. {\it Communications in Statistics- Theory \& Methods}, {\bf43(9)}, 1939-1948.
\bibitem{pn} Psarrakos, G. and Navarro, J. (2013), Generalized cumulative residual entropy
and record values. {\it Metrika}, {\bf76(5)}, 623-640.
\bibitem{ranr} Rajesh, G., Abdul-Sathar, E.I., Nair, K.R.M. and Reshmi, K.V. (2014), Bivariate extension of dynamic cumulative
residual entropy. {\it Statistical Methodology}, {\bf16}, 72-82.
\bibitem{ranr1} Rajesh, G., Abdul-Sathar, E.I., Reshmi, K.V. and Nair, K.R.M. (2014), Bivariate generalized cumulative
residual entropy. {\it Sankhy$\bar{a}$: A}, {\bf76(1)}, 101-122.
\bibitem{r} Rao, M. (2005), More on a new concept of entropy and
information. {\it Journal of Theoretical Probability}, {\bf
18(14)}, 967-981.
\bibitem{rcv} Rao, M., Chen, Y., Vemuri, B.C. and Wang, F.
(2004), Cumulative residual entropy: a new measure of information.
{\it IEEE Transactions on Information Theory}, {\bf 50(6)},
1220-1228.
\bibitem{r} Roy, D. (2002), A characterization of model approach for generating bivariate life distributions using reversed hazard rates. {\it Journal of Japan Statistical Society}, {\bf32(2)}, 239-245.
\bibitem{rcv} Thapliyal, R., Kumar, V. and Taneja, H.C. (2013), On dynamic cumulative entropy of order statistics. {\it Journal of Statistics Applications and Probability}, {\bf2(1)}, 41-46.
\bibitem{sp} Seshadri, V. and Patil, G.P. (1964), A characterization of bivariate distributions by the marginal and conditional distributions
of the same component, {\it Annals of the Institute of Statistical Mathematics}, {\bf15(1)}, 215-221.
\bibitem{ss} Shaked, M. and Shanthikumar, J.G. (2007), {\it Stochastic
Orders}. Springer.
\bibitem{s} Shannon, C.E. (1948), A mathematical theory of
communications. {\it Bell System Technical Journal}, {\bf27},
379-423, 623-656.
\bibitem{sl} Sunoj, S.M. and Linu, M.N. (2012), Dynamic cumulative residual Renyi's entropy. {\it Statistics}, {\bf46(1)}, 41-56.
}
\end{thebibliography}
\end{document}